\documentclass[12pt]{article}
\usepackage{latexsym,color,amsmath,amsthm,amssymb,amscd,amsfonts,epsfig}

\newtheorem{lemma}{Lemma}[section]

\newtheorem{theorem}[lemma]{Theorem}
\newtheorem{definition}[lemma]{Definition}
\newtheorem{corollary}[lemma]{Corollary}

\newtheorem*{remark*}{Remark}

\newcommand{\Id}{{{\mathchoice {\rm 1\mskip-4mu l} {\rm 1\mskip-4mu l}
      {\rm 1\mskip-4.5mu l} {\rm 1\mskip-5mu l}}}}

\begin{document}
\title {  The geometry of partial order on contact transformations of prequantization manifolds    }\author{Gabi Ben Simon
\footnote{This paper is a part of the author's Ph.D. thesis, being
carried out under the supervision of Professor Leonid Polterovich,
at Tel-Aviv University.} \footnote{Partially supported by Israel
Science Foundation grant 11/03.}
\\
School of Mathematical Sciences \\
Tel Aviv University
69978 Tel Aviv, Israel \\
 gabiben@post.tau.ac.il} \maketitle
\begin{abstract} \noindent
In this paper we find  connection between the Hofer's metric of
the group of Hamiltonian diffeomorphisms of a
 closed symplectic manifold, with an integral symplectic form, and the
 geometry, defined in~\cite{El-P}, of the quantomorphisms group
 of its prequantization manifold. This gives two main results: First, we calculate, partly, the geometry
 of the quantomorphisms groups of a prequantization manifolds
 of an integral symplectic manifold which admits certain Lagrangian
 foliation. Second, for every prequantization manifold we give a formula
 for the distance between a point and a distinguished curve in the metric
 space associated to its group of quantomorphisms.
  Moreover, our first result is a
 full computation of the geometry related to the symplectic linear group which can be considered
 as a subgroup of the contactomorphisms group of suitable prequantization manifolds of the
 complex projective space.
 In the course of the proof we use in an essential way the Maslov quasimorphism.


\end{abstract}

\section{Introduction and results}
 The motivating background of this paper can be described  as
follows. Let $(M,\omega)$ be a closed symplectic manifold. Denote
by $Ham(M,\omega)$ the group of all Hamiltonian symplectomorphisms
of $(M,\omega).$ This important group carries a natural
bi-invariant metric called the Hofer metric (a detailed
description of the
 investigation of this metric structure can be found in~\cite{P}).
An analogue of $Ham(M,\omega)$ in the case of contact geometry
 is the group of Hamiltonian
contactomorphisms $Cont(P,\xi)$ where $(P,\xi)$ is a contact
manifold (see~\cite{A} appendix-4,~\cite{AG} and~\cite{AM} for
basic information on contact manifolds). The groups $Ham$ and
$Cont$ are closely related. However in contrast to the case of
$Ham$, no interesting bi-invariant metrics on the group $Cont$ are
known. In~\cite{El-P} Eliashberg and Polterovich noticed that for
certain contact manifold the universal cover of $Cont$ carries a
bi-invariant partial order from which we get a natural metric
space $Z$ associated to $Cont$. The definition of $Z$ is somewhat
indirect, and its geometry is far from being understood even in
the simplest examples.

\noindent In this paper we study the geometry  of $Z$ for certain
subgroups of $Cont(P,\xi)$ where $(P,\xi)$ is a prequanization
space of a closed symplectic manifold $(M,\omega)$. As we remarked
above an important tool is the establishment of a connection
between the geometry of $Z$ and the geometry of the universal
cover of $Ham(M,\omega)$ endowed with the Hofer metric.

\noindent Finally let us mention that further developments on this
subject appear in~\cite{E-K-P}, where in this work we refer mostly
to~\cite{El-P}.

\subsection{Preliminaries on partially ordered groups} \label{basics on p.o.g}
A basic idea in this work is the implementation of basic notions
from the theory of partially ordered groups, to the universal
cover of the groups of contact transformations, of the relevant
contact manifolds. So we start with the following basic
definitions and constructions.

\begin{definition}{\rm Let $\mathcal{D}$ be a group. A subset $\mathcal{C} \subset \mathcal{D}$ is called a\textit{ normal
cone} if

$ \textbf{a.}\text{ } f,g\in \mathcal{C}\Rightarrow fg \in
\mathcal{C}$

$ \textbf{b.}\text{ }f \in \mathcal{C}, h \in \mathcal{D}
\Rightarrow hfh^{-1} \in \mathcal{C}$

$ \textbf{c.}\text{ } \textbf{1}_{\mathcal{D}} \in \mathcal{C}$}

\end{definition}

We define for $f,g \in \mathcal{D}$ that $f \geq g$ if $fg^{-1}\in
\mathcal{C}$. It is not hard to check that this relation is
reflexive and transitive.

\begin{definition}{\rm If the above relation is also anti-symmetric,
then we call it a \textit{ bi-invariant partial order} induced by
$\mathcal{C}$. In this situation, an element $f\in
\mathcal{C}\setminus \mathbf{1}$ is called a \textit{dominant} if
for every $g \in \mathcal{D}$ there exists $n \in \mathbb{N}$ such
that $f^{n}\geq g$.}
\end{definition}

\noindent $\mathbf{Remark.}$ Notice that the normality of the cone
$\mathcal{C}$ implies that for every $f,g,d,e \in \mathcal{D}$
\begin{equation}\text{ if  } f\geq g \text{
 and  } d \geq e \text{  then  } fd \geq ge.
\end{equation}
\begin{definition}{\rm Let $f$ be a dominant and $g\in \mathcal{D}$.
Then the \textit{relative growth} of $f$ with respect to $g$ is
$$\gamma(f,g)=\lim \limits_{n\rightarrow \infty}
\frac{\gamma_{n}(f,g)}{n}$$ where $$\gamma_{n}(f,g)=inf \{p\in
\mathbb{Z}|\text{ } f^{p}\geq g^{n}\}.$$}
\end{definition} The above limit exists as the reader can check by
himself (see also~\cite{El-P} section 1).

Now we want to relate a geometrical structure to the function
$\gamma$ defined above. Denote by $\mathcal{C}^{+}\in\mathcal{C}$
the set of all dominants. We define the metric space $(Z,d)$ in
the following way. First note that $$f,g,h\in
\mathcal{C}^{+}\Rightarrow \gamma(f,h)\leq \gamma(f,g)
\cdot\gamma(g,h).$$ We define the function
$$K:\mathcal{C}^{+}\times \mathcal{C}^{+}\rightarrow [0,\infty)$$
by $$K(f,g)= \max\{\log\gamma(f,g),\log\gamma(g,f)\}.$$ It is
straightforward to check that $K$ is non negative, symmetric,
vanishes on the diagonal and satisfies the triangle inequality.
Thus $K$ is a pseudo-distance.

\noindent Define an equivalence relation on $\mathcal{C}^{+}$ by
setting $f\sim g$ provided $K(f,g)=0$. Put
$$Z=\mathcal{C}^{+}/\sim.$$

The function $K$ on $\mathcal{C}^{+}$ projects in a natural way to
a genuine metric $d$ on $Z$. Thus we get a metric space naturally
associated to a partially ordered groups (see also~\cite{El-P} for
more information).
\subsection{Geometry of the symplectic linear group}\label{linear
geometry}

Let $Sp(2n,{\mathbb{R}})$ be the symplectic linear group. We
denote by $\mathcal{S}$ its universal cover with base point $\Id$,
the identity matrix. We can think of $\mathcal{S}$ as the space of
paths starting at $\Id$ up to a homotopy relation between paths
with the same end point. Throughout this subsection we consider
the group $\mathcal{D}$ to be $\mathcal{S}$.

\noindent Next, consider the equation:\begin{equation}\label{real
positive} \dot{X}(t){X^{-1}(t)}=JH(x,t)\end{equation} Here, $X(t)
\in {\mathcal{S}}$, $H$ is a time dependant symmetric matrix on
$\mathbb{R}^{2n}$ and $J$ is the matrix
$\left( \! \! \! \begin{array}{rr} 0 & -{\Id} \\
 {\Id} &  0 \end{array} \! \! \! \right) $ on $\mathbb{R}^{2n}$.
The quadratic form represented by the matrix $H$, will be called
the {\it Hamiltonian} generating $X(t)$.
 It is easy to verify the following two facts:\\

\noindent {\bf A}). The set of elements $X(t)$ in $\mathcal{S}$,
generated by $H(x,t)$ which are non negative as a quadratic form
for each $t$, establishes a normal cone.

\noindent {\bf B}). In particular, those elements which are
generated by a {\bf strictly} positive $H(x,t)$ are dominants in
$\mathcal{S}$.

\noindent G.I. Olshanskii has proved the following theorem
(see~\cite{Ol} as a reference and~\cite{CSM} as a background
source).

\begin{theorem}\label{linear n.t.p.o}
In the above setting the cone establishes a non trivial partial
order on $\mathcal{S}$.
\end{theorem}In view of theorem~\ref{linear n.t.p.o} we can explore the metric
space $Z$. And indeed we have the following result.

\begin{theorem}\label{main linear theorem} The metric space $(Z,d)$,
derived from the partial order above, is isometric to
$\mathbb{R}$ with the standard metric.
\end{theorem}

\textbf{Remark.}  One can prove theorem~\ref{linear n.t.p.o} using
a criterion developed in~\cite{El-P}. Briefly this criterion says
that in order to establish a non trivial partial order it is
enough to prove that there is no contractible loop in $Sp$ which
is generated by a strictly positive Hamiltonian. We should remark
that the criterion appears in a larger context in~\cite{ El-P}. In
our case the criterion follows from the fact that the Maslov index
of any positive loop is positive, and hence the loop can not be
contractible (see~\cite{M-S} and~\cite{R-S} for more information
on the Maslov index).

\subsection{The geometry of quantomorphisms and Hofer's
metric}\label{pre. hofer}We start with some preliminary
definitions and constructions before presenting the main results
of this section (see~\cite{E-K-P} and~\cite{El-P}).

Let $(P,\xi)$ be a prequantization space of a symplectic manifold
$(M,\omega)$ where $[\omega] \in H^{2}(M,\mathbb{Z})$.
Topologically, $P$ is a principal $S^1$-bundle over $M$. It
carries a distinguished $S^1$-invariant contact form $\alpha$
whose differential coincides with the lift of $\omega$. The
subgroup $Q \subset Cont(P,\xi)$ consisting of all
contactomorphisms which preserve $\alpha$ is called {\it the group
of quantomorphisms} of the prequantization space, (we will give
some more basic information at the beginning of Section~\ref{quant
section}).

We Denote by $\widetilde{Q}$ the universal cover of $Q$, and by
$\widetilde{Ham}(M,\omega)$ the universal cover of $Ham(M,\omega)$
relative to the identity. We have the following important lemma.

\begin{lemma}\label{cw} There is an inclusion \begin{equation}\widetilde{Ham}(M,\omega)\hookrightarrow \widetilde{Q}\end{equation}
of $\widetilde{Ham}(M,\omega)$ into $\widetilde{Q}.$ \end{lemma}

\begin{proof}

Let $\alpha$ be the contact form of the prequantization. The field
of hyperplans $\xi:=\{ker \alpha \}$ is called the contact
structure of $P$. Every smooth function $H:P \times S^{1}
\rightarrow \mathbb{R}$ gives rise to a smooth isotopy of
diffeomorphisms, which preserve the contact structure, as we now
explain. First, define the Reeb vector field of $\alpha$, $Y$, to
be the unique vector field which satisfies the following equations
$$\iota(Y)d \alpha=0 \text{   }, \alpha(Y)=1 .$$

Now given a smooth function $H$ as above define $X_{H}$ to be the
unique vector field which satisfies the following two equations
\begin{equation}\label{good}(1)\text{  }\iota(X_{H}) \alpha =H \text{   }, (2) \text{  }\iota(X_{H})d
\alpha=-dH+(\iota(Y)dH) \alpha \end{equation}

It is well known that the elements of the flow generated by
$X_{H}$ preserve the contact structure. Another important fact is
that the vector field which correspond to the constant Hamiltonian
$H\equiv 1$ is the Reeb vector field, and the flow is the obvious
$S^{1}$ action on $P$.

Now, let $F:M{\times}S^{1}{\rightarrow}\mathbb{R}$ be any
Hamiltonian on $M$ and let $X_{F}$ be its Hamiltonian vector
field. Let $p^{*}F$ be its lift to $P$, by definition $p^{*}F$ is
constant along the fibers of $P$.

We claim that
$$p_{*}X_{p^{*}F}=X_{F}.$$

To see this we substitute $X_{p^{*}F}$ in equation 2
of~\eqref{good}. We get $$\iota(X_{p^{*}F})d
\alpha=-dp^{*}F+(\iota(Y)dp^{*}F) \alpha$$ remembering that $d
\alpha=p^{*}\omega$ we have the equation

$$\iota(X_{p^{*}F})p^{*}\omega=-dp^{*}F+(\iota(Y)dp^{*}F) \alpha$$

$\Longleftrightarrow$

For every tangent vector to $P$, $v$, we have

$$\iota(p_{*}X_{p^{*}F})\omega(p_{*}v)=-dp^{*}F(v)+(\iota(Y)dp^{*}F) \alpha(v)$$

$\Longleftrightarrow$

\begin{equation}\label{almost}\iota(p_{*}X_{p^{*}F})\omega(p_{*}v)=-p^{*}dF(v)+(\iota(Y)dp^{*}F)
\alpha(v).\end{equation}

Now note that $\iota(Y)dp^{*}F\equiv 0.$ This is since $Y$ is
tangent to the fibers of $P$, and since $p^{*}F$ is constant along
the fibers, its differential has no vertical component with
respect to any local system of coordinates.

Thus, we get from equation~\eqref{almost} the equation

$$\iota(p_{*}X_{p^{*}F})\omega(p_{*}v)=-dF(p_{*}v)$$
for every tangent vector $v$. Now since the linear map $p_{*}$ is
surjective on every point of $P$, we get the equation
$$\iota(p_{*}X_{p^{*}F})\omega=-dF.$$
The fact that $\omega$ is non degenerate gives the desired
equation $p_{*}X_{p^{*}F}=X_{F}.$

 Denote by $\mathcal{F}$ the set of all time-1-periodic functions

 \begin{equation}F:M{\times}S^{1}{\rightarrow}\mathbb{R}\text{    such that   }
 \int_{M}{F(x,t)\omega^{n}}=0 \text{ for every  } t \in
 [0,1].\end{equation}

 It is not hard to check that for every representative,
$\{f_{t}\}_{t\in[0,1]},$ of an element of
$\widetilde{Ham}(M,\omega)$ there is a unique Hamiltonian
$F\in\mathcal{F}$ which generates it.

 Let $F$ be any Hamiltonian. The lift of $F$ is the contact Hamiltonian
$$\widetilde{F}:P{\times}S^{1}{\rightarrow}\mathbb{R}$$ defined by $\widetilde{F}:=p^{*}F$
where $p:P \rightarrow M$ is the projection from $P$ to $M$  (Note
that the function $\widetilde{F}$ is constant along the fibers of
$P$).  Now, $\widetilde{F}$ is a contact Hamiltonian function of
the contact manifold $P$ which generates an element of
$\widetilde{Q}$. It is well known that every representative of an
element of $\widetilde{Q}$, is generated by a unique contact
Hamiltonian which is a lift of an Hamiltonian defined on $M.$ We
further remind the reader that $Q$ is a central extension of
$Ham$. That is we have the short exact sequence $$\Id \rightarrow
S^{1} \rightarrow Q \rightarrow \text{  } Ham \text{  }
\rightarrow  \Id .$$

Which is from the Lie algebra point of view, the Poisson bracket
extension of the algebra of Hamiltonian vector fields. That is

$$0 \rightarrow \mathbb{R} \rightarrow  C^{\infty}(M) \rightarrow
\text{  }\mathcal{F} \text{  } \rightarrow 0 .$$

We first explain how to map the identity element of
$\widetilde{Ham}$, $[\Id]$, into $\widetilde{Q}$. Let $f_{t}$ be
any contractible loop representing $[\Id]$, and let $F_{t}$ be the
unique normalized Hamiltonian which generate it. The path
$\hat{f}_{t}$, in $Q$, which is generated by the lift of $F_{t}$,
$p^{*}F_{t}$, is a lift of $f_{t}$. Nevertheless, it is not
necessarily closed (actually it is closed, as we show below). It
easy to check that there is a constant, $c$, such that
$p^{*}F_{t}+c$ generates a lift, $\hat{h}(t)$, of the loop $f_{t}$
and this lift is a loop based at the identity of $Q.$ We claim
that

\textbf{a.} The loop $\hat{h}(t)$ is contractible.

\textbf{b.} The constant $c$ equals zero.

We first prove $a$. Since the loop $f_{t}$ is contractible, there
is a homotopy $$K_{s}(t):[0,1] \times [0,1] \rightarrow Ham$$ such
that for every $s$, $K_{s}(t)$ is a contractible loop. We denote
by $F_{s}(x,t)$ the normalized Hamilonian which generate the loop
$K_{s}$. As we explained, for every $s$ there is a constant $c(s)$
such that $p^{*}F_{s}+c(s)$ generates a loop, $\hat{h}_{s}(t)$, in
$Q$ which is a lift of $K_{s}.$

Now since  $\hat{h}_{s}(t)$ is a lift of the homotopy $K_{s}(t)$,
(which is a homotopy between  $f_{t}$ and $\Id$) it is a homotopy
between $\hat{h}(t)$ and the identity element of $Q.$ Thus
$\hat{h}(t)$ is contractible.

We now prove $b$. First we remind the reader the Calabi-Weinstein
invariant, which is well defined on the group $\pi_{1}(Q)$. It was
proposed by Weinstein in~\cite{W1}. For a loop $\gamma \in Q$
define the Calabi Weinstein invariant by
$$cw(\gamma)=\int_{0}^{1}dt \int \limits_{M}F_{t} \omega^{n}$$

where $p^{*}F_{t}$ is the unique Hamiltonian generating $\gamma.$

Now it is clear that $$cw([\Id])=0.$$

We conclude that for every $s$ $$\int_{0}^{1}dt \int
\limits_{M}(F_{s}(t,x)+c(s)) \omega^{n}=0$$

$\Longleftrightarrow$
$$c(s)vol(M)=-\int_{0}^{1}dt \int
\limits_{M}F_{s}(t,x) \omega^{n}.$$

Now since $F_{s}$ is normalized for every $s$ we get that $c(s)=0$
as desired.

The map of an arbitrary element of $\widetilde{Ham}$ is defined in
the same manner. Let $f_{t}$ be any representative of an element,
$[f]$, of $\widetilde{Ham}.$ Let $F_{t}$ be the unique normalized
Hamiltonian generating $f_{t}$. Then the image of $[f]$ is the
homotopy class of the path, $\hat{f}_{t}$, generated by the lift
of the Hamiltonian $F_{t}$, namely $p^{*}F_{t}$. By using the same
type of argument as above it can be shown that this map is well
defined. Finally, this map is clearly injective.

\end{proof}

The universal cover of $Q$ carries a natural normal cone $C$
consisting of all elements generated by non-negative contact
Hamiltonians. The set of dominants consists of those elements of
$\widetilde{Q}$ which are generated by strictly positive
Hamiltonian. The normal cone always gives rise to a genuine
partial order on $\widetilde{Q}$. Thus one can define the
corresponding metric space $(Z,d)$ as it was explained in
subsection \ref{basics on p.o.g}. We state this as a theorem.

\begin{theorem} Let $(M,\omega)$ be a symplectic form such that
$\omega \in H^{2}(M,\mathbb{Z})$. Let $(P,\xi)$ be a
prequantization of $(M,\omega).$ Then $\widetilde{Q}$, the
universal cover of the group of quantomorphisms of $Q$, is
orderable.

\end{theorem}

\begin{proof} According to ~\cite{El-P} all we need to check is
that there are no contractible loops in $Q$ which are generated by
a strictly positive contact Hamiltonian. Such contractible loops
can not exists since the Calabi-Weinstein invariant of any
representative of a contractible loop must be zero, and thus can
not be generated by a strictly positive Hamiltonian.

\end{proof}

In light of theorem~\ref{linear n.t.p.o} we should remark the
following.

 \noindent \textbf{Remark.}  Let $(M,\omega)$ be the
 complex projective space $\mathbb{C}P^{n-1}$ endowed with
 the Fubini-Study symplectic form normalized to be integral. Let $(P,\xi)$ be the sphere
 $S^{2n-1}$ endowed with the standard contact structure, and let $(\mathbb{R}P^{2n-1},\beta)$
 be the standard contact real projective space. It is well known that these two contact manifolds are
 prequantizations
 of $(M,\omega)$ (see~\cite{Ki},~\cite{W} and~\cite{Wo} for preliminaries on
 prequantization). Let
   $Cont(P,\xi)$ and $Cont(\mathbb{R}P^{2n-1},\beta)$ be the groups of all Hamiltonian contact transformations of
  these manifolds.

\noindent The relations of the symplectic linear group to these
groups is as follows. The group $ Sp(2n,\mathbb{R})$ is  subgroup
of all elements of $Cont(P,\xi)$ which commutes with $-\Id$, and $
Sp(2n,\mathbb{R})/\pm \Id$ is a subgroup of
$Cont(\mathbb{R}P^{2n-1},\beta).$
  Now, it is known that $Cont(\mathbb{R}P^{2n-1},\beta)$ admits a nontrivial partial order
  (while $Cont(P,\xi)$
  does not). The non orderabilty of  $Cont(P,\xi)$
  follows from~\cite{E-K-P}, this is done by
  using the criterion mentioned above. The orderabilty of $Cont(\mathbb{R}P^{2n-1},\beta)$ follows from the theory of the
  nonlinear Maslov index introduce by Givental see~\cite{G}.

\bigskip In what follows we establish a connection between the
partial order on $\widetilde{Q}$ (and its metric space $Z$) and
the Hofer's metric. So first let us recall some basic definitions
related to the Hofer's metric. The \textit{Hofer's distance}
between  $f \in \widetilde{Ham}(M,\omega)$ and $\Id$ is defined by
$$\rho(f,\Id)=\inf_{G} \int_0^1 \left ( \max_{x{\in}M}G(x,t)-\min_{x{\in}M}G(x,t) \right )dt, $$
where the infimum is taken over all time $1$-periodic Hamiltonian
functions $G$ generating paths in $Ham$ which belongs to the
homotopy class represented by $f$ with fixed end points
 (see~\cite{H-Z},~\cite{M-S} and~\cite{P}, for further
information). Furthermore, we define the\textit{ positive and the
negative part of the Hofer's distance} as
$$\rho_{+}(f,\Id) = \inf { \int_0^1 \max_{x \in M} {G(x,t)}}dt, \ \ \ \ \rho_{-}(f,\Id)=
\inf{\int_0^1-\min_{x{\in}M}{G(x,t)}dt},$$ 
where the infimum is taken as above.

 \noindent Now, for each $f \in \widetilde{Ham}(M,\omega)$, generated by some Hamiltonian, set
 \begin{equation} \| f \|_{+} = \inf \left \{ \max F(x,t) \right \} , \ \ \ \| f \|_{-} =
 \inf \left \{ - \min F(x,t) \right \}\end{equation}
 where the infimum is taken over all Hamiltonian functions $F \in
 \mathcal{F}$ generating paths in $Ham$ which belongs to the
homotopy class represented by $f$.
   The following result is due to Polterovich
see~\cite{P-1}.

 \begin{lemma} \textit{For every $f \in
\widetilde{Ham}(M,\omega)$ we have
\begin{equation}\label{+ equation}\! \! \! \rho_{+}(f)=\|f\|_{+}\end{equation}
\begin{equation}\rho_{-}(f)=\|f\|_{-}\ . \end{equation}}
\end{lemma}

We further define the positive and negative \textit{asymptotic}
parts of the Hofer's metric:
$$\|f\|_{+,\infty}:=\lim\limits_{n\rightarrow
\infty}\frac{\|f^{n}\|_{+}}{n}\text{ and }
\|f\|_{-,\infty}:=\lim\limits_{n\rightarrow
\infty}\frac{\|f^{n}\|_{-}}{n}.$$

Before stating our first main result we fix the following
notation.  Given $(P,\xi)$  (a prequantization of a symplectic
manifold $(M,\omega)$) we denote by $e^{is}$ the diffeomorphism of
$P$ obtained by rotating the fibers in total angle $s.$
 Note that the set $\{e^{is}\}_{s\in \mathbb{R}}$ is a one parameter
family of contact transformations in $Q.$

\begin{theorem}\label{main theorem 1}

 Let $f$ be the
time-1-map of a flow generated by  Hamiltonian $F \in
\mathcal{F}$. Let $\widetilde{F}$ be the lift of $F$ to the
prequantization space and let $\tilde{f}$ be its time-1- map. Take
any $s \geq 0$ such that $e^{is}\tilde{f}$ is a dominant. Then
\begin{equation}\label{log}dist(\{e^{it}\},e^{is}\tilde{f}):=\inf\limits_{t}K(e^{it},e^{is}\tilde{f})=
\frac{1}{2}\log\frac{s+\|f\|_{+,\infty}}{s-\|f\|_{-,\infty}}.\end{equation}
\end{theorem}

Before stating our second result we need the following
preliminaries.

\begin{definition}{\rm  Let $(M,\omega)$ be a closed symplectic manifold.
Let $L$ be a closed Lagrangian submanifold of $M$. We say that $L$
has the \textit{Lagrangian intersection property} if $L$
intersects its image under any exact Lagrangian isotopy.}
\end{definition}

\begin{definition}{\rm  We say that \textit{L has the stable Lagrangian
intersection property} if $L\times \{r=0\}$ has the Lagrangian
intersection property in $(M\times T^{\ast}S^{1},\omega \oplus
dr\wedge dt)$ where $(r,t)$ are the standard coordinates on the
symplectic manifold $T^{\ast}S^{1}$ and $M\times T^{\ast}S^{1}$ is
considered as a symplectic manifold with the symplectic form
$\omega \oplus dr\wedge dt$.}
\end{definition}

For more information on Lagrangian intersections see~\cite{P}
chapter 6.

Consider now the following situation. For $(M,\omega)$ a closed
symplectic manifold, $[\omega] \in H^{2}(M,\mathbb{Z})$, assume
that we have an open dense subset $M_{0}$ of $M$, such that
$M_{0}$ is foliated by a family of closed Lagrangian submanifolds
$\{L_{\alpha}\}_{\alpha\in \Lambda}$, and each Lagrangian in this
family has the stable Lagrangian intersection property.

We denote by \texttt{F} the family of autonomous functions, in
$\mathcal{F}$, defined on $M$ which are constant when restricted
to $L_{\alpha}$ for every $\alpha,$ note that this means that the
elements of \texttt{F} commutes relative to the Poisson brackets.
That is for every $F,G\in \texttt{F}$ we have $\{F,G\}=0.$

 We denote by $\{\texttt{F},\|\text{ }\|_{max}\}$ the metric space of the family
of functions \texttt{F} endowed with the $max$ norm where the
$max$ norm is $\|F\|_{max}:=\max\limits_{x\in M}|F(x)|.$

In the course of the proof we will use the following subspace: Let
$V$ be the subspace of $\widetilde{Q}$ constitutes of all elements
generated by contact Hamiltonians of the form $s+\widetilde{F}$
where $\widetilde{F}$ is a lift of an Hamiltonian $F\in
\texttt{F}$ and $s+\widetilde{F}> 0.$

We now state our second result.

 \begin{theorem}\label{main theorem}

   There is an isometric injection of the metric space
$(\texttt{F},\|\text{ }\|_{max})$ into the metric space $Z$.

\end{theorem}

Note that the existence of $s$ in theorem~\ref{main theorem} is
justified due to the fact that the function $F$, and thus
$\widetilde{F}$, attains a minimum.

As we will show in section~\ref{quant section} theorem~\ref{main
theorem} can be implemented to the case of the standard even
dimensional symplectic torus
$T^{2n}=\mathbb{R}^{2n}/\mathbb{Z}^{2n}$. Our second example will
be for oriented surfaces of genus greater equal 2. In this case we
foliate the surface (minus a finite family of curves) with closed
non contractible loops which are all closed lagrangians with the
stable Lagrangian intersection property.

\section{ Proof of theorem~\ref{main linear theorem}}
\subsection{The unitary case}\label{unitary case}
In this subsection we prove theorem~\ref{main linear theorem} for
the unitary group $U(n)$. This is a result by itself. Moreover
elements of the proof will be used to prove theorem~\ref{main
linear theorem} for the symplectic case. Here we use occasionally
basic properties of the exponential map on the Lie algebra of a
matrix Lie group (see~\cite{H} chapter 2).

We denote by $\mathcal{A}$ the universal cover of $U(n)$. We view
$\mathcal{A}$ as the set of all paths starting at $\Id$, the
identity element, up to homotopy relation between paths with the
same end point. We will denote an element $[u]\in \mathcal{A}$,
with a slight abuse of notation, by $u$ where $u$ is a path
representing $[u]$.

As in the symplectic case (see definition~\ref{maslov quasi}) we
can define for an element $u \in \mathcal{A}$ its Maslov index as
$\alpha(1)-\alpha(0)$ where $e^{i2 \pi \alpha(t)}$= det$(u(t))$.
We will denote it by $\mu(u)$. Using the same definitions
of~\ref{linear geometry} we have the following. An element $u \in
\mathcal{A}$ is called \textit{semi positive} if for some
representative, the hermitian matrix $h_{u}$ defined by the
equation
\begin{equation}\label{complex positive} \dot{u}u^{-1}=ih_{u} \end{equation}
is semi positive definite. Note that if we consider $U(n)$ as a
subgroup of $Sp(2n,\mathbb{R})$ via the
realization$$ A+iB\mapsto \left( \! \! \! \begin{array}{rr} A & {-B} \\
 {B} &  A \end{array} \right)  $$
 equation~\eqref{complex positive} is the same as equation~\eqref{real
 positive}. We proceed as in~\ref{linear geometry}, define a
 partial ordering on $\mathcal{A}$, $\geq$, by $f \geq g$ if and
 only if $fg^{-1}$ is semi positive. Let $\mathcal{C}$ denote the
 subset of all semi positive elements and denote by $\mathcal{C^{+}}$ the set of positive
 definite elements
  of $\mathcal{A}$. These are, respectively, the normal
 cone and the set of dominants  of $\mathcal{A}$ ( see~\ref{basics on p.o.g}).

  We define the function $\gamma_{n}$, the relative growth, and the metric
 space $Z$ in the same way as in~\ref{basics on p.o.g}. 

\subsubsection{Preliminary basic lemmas}

 We need
the following two elementary lemmas.
\begin{lemma}\label{homotopy}Let $u$ and $v$ be two paths in $U(n)$ with the same
endpoints. Then $[u]=[v]$ if and only if $[|u|]=[|v|]$ in $U(1)$.

\end{lemma}

\begin{proof} It is well known that the determinant map
$U(n)\rightarrow U(1)$ induces an isomorphism between the
fundamental groups of these spaces. Let $\hat{\gamma}$ be the
inverse path to $\gamma$ and let $\ast$ be the juxtaposition of
paths. Thus $[u]=[v]$ if and only if $u \ast \hat{v}$ is homotopic
to the identity, if and only if $$|u \ast \hat{v}| = |u| \ast
\widehat{|v|}$$ is homotopic to the identity, if and only if
$[|u|]=[|v|].$
\end{proof}

\begin{lemma}\label{tr lemma} Assume that $u$, $h_{u}$ satisfies
equation~\eqref{complex positive}. Then \begin{equation}
\mu(u)=\int_{0}^{1}\text{tr}h_{u}(t)dt .\end{equation}

\end{lemma}

\begin{proof} For a complex path $z(t)\neq0$ we have of course
$\frac{d}{dt}\log z(t)=\frac{\dot{z}(t)}{z(t)}$. Taking the
imaginary part we have:
$$\frac{d}{dt}argz(t)=Im\frac{\dot{z}(t)}{z(t)}.$$ Fix $t=t_{0}$
in $[0,1]$ and write $u(t)=u(t_{0})v(t)$ where $v(t_{0})=\Id$. Let
$D$ be the differentiation with respect to $t$ at $t=t_{0}$. Note
that $$\dot{v}v^{-1}=ih_{u}.$$ Then $D|v|=\text{tr}(\dot{v})$ at
$t_{0}$, which implies that $$D\text{ }
arg|u|=Im\frac{D|u|}{|u|}=Im\frac{D|v|}{|v|}$$ $$=Im\text{
}\text{tr}( \dot{v})=Im\text{ }\text{tr}(ih_{u}v)=\text{tr}
(h_{u}).$$

In the computation above remember that $v=1$ at $t=t_{0}$. Now the
above formula follows as $arg|u(0)|=0.$

\end{proof}

As a corollary of the lemma we have that if $u$ is semi positive
then $\mu(u)\geq0$ (since if $h_{u}$ is semi positive then
$\text{tr}h_{u} \geq 0$). 
The converse of the corollary is not true (we will not give an
example). Nevertheless we have:

\begin{lemma}\label{main linear lemma} If $\mu(v) \geq 2\pi n$ then $v \geq \Id$. 
In this case, we can have a representative
$v(t)=e^{itA}$ where $\text{tr}(A)=\mu(v)$, $A$ is hermitian and
positive semidefinite, and if $\lambda_{1},...,\lambda_{n}$ are
the eigenvalues of $A$, then
$\max_{i,j}|\lambda_{i}-\lambda_{j}|\leq2 \pi n.$

\end{lemma}

\begin{proof} Suppose that $\mu(v) \geq 2\pi n$. Then according to
lemmas~\ref{homotopy} and~\ref{tr lemma} we can write
$v(t)=e^{itA}$ with $A$ hermitian. The choice of $A$ is not
unique, however $\text{tr}(A)=\mu(v)$ up to a multiple of $2\pi$.
Still, we can modify $A$ so that $\text{tr}(A)=\mu(v)$. To see
this, diagonalize $A$ by a unitary matrix and then modify the
eigenvalues by multiples of $2 \pi$ (here we are using the fact
that if $C$ is an invertible matrix then
$e^{CXC^{-1}}=Ce^{X}C^{-1}$ for any arbitrary matrix $X$).

In fact we can do better: We may still modify without changing
$\text{tr}A$ such that now $A$ will be positive semidefinite and
$\max_{i,j}|\lambda_{i}-\lambda_{j}|\leq2 \pi n.$ To see this,
after $A$ has been diagonalize, let $\lambda_{1},...,\lambda_{n}$
be the eigenvalues of $A$ (in an increasing order). Then we have
$$\sum \limits_{i}\lambda_{i} = \text{tr}A= \mu(v) \geq 2\pi
n.$$

\noindent Now let $0< \bar{\lambda}_{i}\leq 2 \pi$ be the unique
number that is congruent to $\lambda_{i}$ modulo $2 \pi
\mathbb{Z}.$  Then we have

$$\sum \limits_{i}\bar{\lambda}_{i}=\text{tr}A-2\pi k,\text{
}k>0.$$ Note that
$\max_{i,j}|\bar{\lambda}_{i}-\bar{\lambda}_{j}|\leq2 \pi .$ We
first assume that $k<n$ (the case $n=1$ is trivial so we assume $n
\geq 2$). In this case we change the $\bar{\lambda}_{i}\text{'s}$
by distributing the $k$ $2 \pi \text{'s}$ to the first $k$
$\bar{\lambda}_{i}\text{'s}.$ It is clear that the new maximal
$\bar{\lambda}_{i}$ will come from the first, new, $k$
$\bar{\lambda}_{i}\text{'s}$ and that the condition
$\max_{i,j}|\bar{\lambda}_{i}-\bar{\lambda}_{j}|\leq2 \pi n$ is
kept.

Now if $k \geq n$ then write $k=l n+m$ where $l,m \in \mathbb{N}$
and $m<n$. In this case we add $2 \pi$ to all the
$\bar{\lambda}_{i}\text{'s}$ $l$ times and then we add to the
first $m$ $\bar{\lambda}_{i}\text{'s}$ the remaining $m$ $2 \pi
\text{'s}.$ Now we are at the same position as in the first case
where $k < n.$ The resulting matrix $A$ (before diagonalization)
has the desired properties.

\end{proof}

\subsubsection{The Metric induced by $\mathcal{C}^{+}$}\label{c^+}

Let $f,g \in {\mathcal{C}}^{+}.$ We wish to compute $\gamma(f,g).$
We have the following theorem.
\begin{theorem}\label{hand theorem} For every non constant $f,g\in
{\mathcal{C}}^{+}$ we have: \begin{equation}
\gamma(f,g)=\frac{\mu(g)}{\mu(f)}.\end{equation}

\end{theorem}
\begin{proof} 

We start with the following definition.
\begin{equation}\label{gamma*}\gamma^{*}(f,g)=\inf\{\frac{r}{s}|f^{r}\geq
g^{s}\text{ , }r \in \mathbb{Z}\text{ , }s \in
\mathbb{N}\}\end{equation} We claim that $\gamma^{*}=\gamma$.
Indeed, denote by $T$ the set of numbers which satisfies the
condition of the right hand side of~\eqref{gamma*}. Assume that $
\frac{r}{s}\in T $, then the equivalence of the definitions
follows from the inequality
$$r\geq \gamma_{s}(f,g)\geq s\gamma^{*}(f,g).$$ We claim that a
sufficient condition for $f^{r}\geq g^{s}$ is that $r
\mu(f)-s\mu(g) \geq 2 \pi n.$ Indeed, if so then
$\mu(f^{r}g^{-s})=r \mu(f)-s\mu(g) \geq 2 \pi n,$ so by
lemma~\ref{main linear lemma}, $$f^{r}g^{-s} \geq
\Id\Leftrightarrow f^{r}\geq g^{s}.$$ Now, given $\varepsilon >0$
choose $s> \frac{1}{\varepsilon}.$ Let $r \geq 0$ be the smallest
integer so that $$r \mu(f)-s\mu(g) \geq 2 \pi n.$$ Thus, $$r
\mu(f)-s\mu(g)<2 \pi n+\mu(f)\Leftrightarrow$$ $$\frac{r}{s}
\mu(f)-\mu(g) <(2 \pi n+\mu(f))\varepsilon \Leftrightarrow$$
$$\frac{r}{s}< \frac{\mu(g)}{\mu(f)}+\frac{(2 \pi
n+\mu(f))\varepsilon}{\mu(f)}.$$ Since $\varepsilon$ is arbitrary
then $$\gamma^{*}(f,g)=\inf\{\frac{r}{s}|f^{r}\geq g^{s}\} \leq
\frac{\mu(g)}{\mu(f)} .$$ Now $$\gamma(f,g) \gamma(g,f) \geq 1$$
gives the desired equality

$$\gamma(f,g)=\frac{\mu(g)}{\mu(f)}.$$

\end{proof}

Now we can finish the proof of theorem~\ref{main linear theorem}
in the unitary case. First, we remind the definition of the metric
$Z$. $Z=\mathcal{C}^{+}/\sim $ where $f\sim g$ provided $K(f,g)=0$
and $K(f,g)= \max\{\log\gamma(f,g),\log\gamma(g,f)\}.$ Define a
map $p:\mathcal{A}\rightarrow \mathbb{R}$ by $p(u)=\log(\mu(u)).$
We claim that the map $p$ induces an isomorphism of metric spaces
that is
$$\frac{\mathcal{C}^{+}}{\sim}=Z \cong \mathbb{R}.$$ Indeed, by
theorem~\ref{hand theorem} we get
$$|p(f)-p(g)|=\max \{\log \mu(g)/\mu(f), \log \mu(f)/\mu(g)\}=K(f,g).$$
Thus $p$ is an isometry.

\noindent \textbf{Remark.} Note also that $p$ preserve order as
well. Indeed, $f\geq g$ implies that $0\leq
\mu(fg^{-1})=\mu(f)-\mu(g)$ which implies that $p(f)\geq p(g).$
See~\cite{El-P} subsection 1.7 for more details on this phenomena,
in a larger context.

\subsection{The Maslov qusimorphism}\label{quasi maslov}

Before proving theorem~\ref{main linear theorem} in the symplectic
case we define an important property which we use in a crucial way
in the course of the proof.

\begin{definition} \label{quasi}
{\rm Let $G$ be a group. A \textit{quasimorphism} $r$ on $G$ is a
function $r:G\rightarrow\mathbb{R}$ which satisfies the
homomorphism equation up to a bounded error: there exists $R>0$
such that
$$|r(fg)-r(f)-r(g)|\leq R$$
for all $f,g\in G$}.
\end{definition} Roughly speaking a quasimorphism is a
homomorphism up to a bounded error. See ~\cite{Ba} for
preliminaries on quasimorphisms. A quasimorphism $r_{h}$ is called
\textit{homogeneous} if $r_{h}(g^{m})=mr_{h}(g)$ for all $g\in G$
and $m\in \mathbb{Z}$. Every quasimorphism $r$ gives rise to a
homogeneous one
$$r_{h}(g)=\lim \limits_{m\rightarrow\infty}\frac{r(g^{m})}{m}.$$

As we said a basic tool in the proof of Theorem \ref{main linear
theorem} is the {\it{Maslov quasimorphism}}  whose definition is a
generalization of the Maslov index from loops to paths in
$Sp(2n,{\mathbb{R}})$.

\begin{definition}\label{maslov quasi}{\rm  Let $\Psi(t)=U(t)P(t)$ be a representative of a point in $\widetilde{Sp}$
,where $U(t)P(t)$ is the polar decomposition in $\widetilde{Sp}$
of $\Psi(t)$, $U(t)$ is unitary and $P(t)$ is symmetric and
positive definite. Choose $\alpha(t)$ such that $e^{i2\pi
\alpha(t)}$= det$(U(t))$. Define the Maslov quasimorphism $\mu$ by
$\mu([\Psi])=\alpha(1)-\alpha(0)$ }.
\end{definition}
\begin{theorem}\label{the quasi maslov} The Maslov quasimorphism is a
quasimorphism.
\end{theorem}
For the proof of theorem~\ref{the quasi maslov} see~\cite{BG}
and~\cite{D}. \noindent We denote by $\widetilde{\mu}$ the
homogeneous quasimorphism corresponding to $\mu$:
$${\widetilde{\mu}}(x)=\lim_{k{\rightarrow}\infty}{\frac{\mu(x^k)}{k}}. $$
We further remark that the restriction of the (homogeneous) Maslov
quasimorphism to  $\mathcal{A}$ (since there is an isomorphism
$\pi_{1}(Sp(2n,\mathbb{R}))\cong \pi_{1}(U(n))$, we may consider
$\mathcal{A}$ as a subset of $\widetilde{Sp}$) is the Maslov index
we have defined in subsection~\ref{unitary case} on $\mathcal{A}$.

\subsection{The Symplectic case}

We now turn to the symplectic case. Let $J$ be the standard
symplectic structure on $\mathbb{R}^{2n}$ (see~\ref{linear
geometry}). We recall that the symplectic linear group
$Sp(2n,\mathbb{R})$, is the group of all matrices $A$ which
satisfies $A^{T}JA=J.$ As we already remarked $U(n)$ is a subgroup
of $Sp(2n,\mathbb{R}).$ The elements of $U(n)$ are precisely those
which commute with $J$.

 \noindent We denote by $\mathcal{S}$ the universal
cover of this group having as base point the identity matrix
${\Id}$. As before we use the same letter to denote an element in
$\mathcal{S}$ and a representing path. However we shall use
capital letters instead.

Now, assume that $X,Y \in \mathcal{S}$. Let $H_{X}, H_{Y}, H_{XY}$
be the Hamiltonians generating respectively $X,Y, XY$
(see~\ref{linear geometry}). For later use we need the following
formula. The formula is well known and it is easy to prove.
\begin{equation}\label{prod. formula} H_{XY}=
H_{X}+{X^{-1}}^{T}H_{Y}X^{-1}. \end{equation}

We denote by $\mathcal{S}^{+}$ the set of dominants of
$\mathcal{S}$. Note that we have the inclusion
$\mathcal{C}^{+}\hookrightarrow \mathcal{S}^{+}.$ Recall that
$\mathcal{S}^{+}$ induces a non trivial partial order on
$\mathcal{S}$ (see~\ref{linear geometry}), thus we can define for
every $f,g \in\mathcal{S}^{+}$ the function $\gamma(f,g)$. Here
one must be cautious as this magnitude may have two different
meanings for unitary $f,g$. According to the following theorem
this causes no difficulties (see the remark following
theorem~\ref{sp lemma}).

Recall that $\tilde{\mu}$ denote the homogenous Maslov
quasimorphism. The key theorem of this subsection is the
following.
\begin{theorem}\label{sp lemma} For all $X,Y \in \mathcal{S}^{+}$ we
have:\begin{equation}\label{main symp.
step}\gamma(X,Y)=\frac{\tilde{\mu}(Y)}{\tilde{\mu}(X)}.\end{equation}

\end{theorem}

 Theorem~\ref{main linear theorem}, now, follows from
 theorem~\ref{sp lemma}. Indeed define $p(X)=\log(\tilde{\mu}(X))$
for
 every $X \in \mathcal{S}^{+}.$ Now repeat verbatim the last part
 of~\ref{c^+}.

\textbf{Remark.} Notice that if $X$ and $Y$ are unitary, the
theorem shows that their symplectic relative growth is the same as
the unitary relative growth, since on unitary paths,
$\tilde{\mu}=\mu.$

We now prove theorem~\ref{sp lemma} in two steps.
\subsubsection{Proof of Theorem~\ref{sp lemma} - Step 1}
In step 1 we prove the following lemma.
\begin{lemma} For every positive definite symmetric symplectic
\textbf{matrix} $P$, there exists a positive path $X$ connecting
$\Id$ with $P$, such that \begin{equation}\mu(X)\leq 4 \pi
n.\end{equation}
\end{lemma}
\begin{proof} Without the loss of generality we may assume that
$P$ is diagonal. Otherwise, diagonalize $P$ by a unitary matrix
$U_{0}$ and notice that for every path $Y$, we have by a direct
calculation that $$H_{U_{0}YU_{0}^{-1}}=U_{0}H_{Y}U_{0}^{-1}$$
where $H_{Y}$ is a Hamiltonian generating $Y$. Thus $X$ is
positive if and only if $U_{0}XU_{0}^{-1}$ is.  Now, assume that
the lemma has been proved for the case $n=1$. Using this
assumption we prove the general case.

We first define $n$ different embeddings
$$j_{i}:Sp(2,\mathbb{R})\hookrightarrow Sp(2n,\mathbb{R})\text{,        }1\leq i \leq n$$
as follows. We copy a $2\times 2$ matrix $A$ to the $(i,n+i)$
block of a $2n\times 2n$ matrix $B$. That is
$$B_{ii}:=A_{11},\text{  }B_{ii+n}:=A_{12}$$ $$B_{i+ni}:=A_{21},\text{  }B_{i+ni+n}:=A_{22}.$$
 The rest of the elements of $B$ are defined as follows (we assume
 of course that $(k,l)\neq (i,i),(i,i+n),(i+n,i),(i+n,i+n)$).
 $$B_{kl}:=\begin{cases}0 &k\neq l \\1&k=l\text{ .}\end{cases}$$

 These are at least homomorphisms to $GL(2n,\mathbb{R})$, but in
 fact one can easily check that the image is symplectic. These
 embeddings preserve transpose operator and thus symmetry and
 orthogonality. As a result the embeddings respect
 symmetric-unitary decomposition. In particular
 \begin{equation}\label{mu}\mu(X)=\mu(j_{i}(X))\end{equation}

 \noindent for every $X\in Sp(2,\mathbb{R})$ and $1\leq i \leq n$.

Note the following property of the embeddings $j_{i}$. For every
$i\neq k$ and any $X,Y\in Sp(2,\mathbb{R})$ we have
\begin{equation}\label{j commu.} j_{i}(X)j_{k}(Y)=j_{k}(Y)j_{i}(X)\end{equation}

 Now, by assumption $P$ is diagonal, let
 $(\lambda_{1},...,\lambda_{n},1/\lambda_{1},...,1/\lambda_{n})$
 be the diagonal of $P$. Define
 $P_{i}=diag(\lambda_{i},1/\lambda_{i})$, and notice that
 \begin{equation}\label{prod.1}P=\prod\limits_{i=1}^{n}j_{i}(P_{i})\end{equation}

 Having proved the theorem for $n=1$, we can find paths $X_{i}$ connecting
$\Id$ (of $Sp(2,\mathbb{R})$) with $P_{i}$, such that
$\mu(X_{i})\leq 4 \pi$. Define
\begin{equation}\label{prod.2}X=\prod\limits_{i=1}^{n}j_{i}(X_{i})\end{equation}
We write $j_{i}(X_{i})=P_{i}(X_{i})U_{i}(X_{i})$ for the polar
decomposition of the matrix $j_{i}(X_{i})$ in $Sp(2n,\mathbb{R})$.
Note that according to what we have remarked above, this
representation of this polar decomposition is the image of the
polar decomposition of $X_{i}$ with respect to the homomorphism
$j_{i}$ for all $i$. Anyhow we get
\begin{equation}\label{polar prod.}X=\prod\limits_{i=1}^{n}P_{i}(X_{i})U_{i}(X_{i})\end{equation}
Note that according to~\eqref{j commu.} all the elements in the
product of~\eqref{polar prod.} commutes. Thus we can write
\begin{equation}\label{decom.polar.}X=\prod\limits_{i=1}^{n}P_{i}(X_{i})\prod\limits_{i=1}^{n}U_{i}(X_{i})\end{equation}
We claim that the r.h.s. of~\eqref{decom.polar.} is the polar
decomposition of $X$. Indeed, all the elements in the product
$\prod\limits_{i=1}^{n}P_{i}(X_{i})$, according to~\eqref{j
commu.}, commutes. Thus, this product is a symmetric (positive
definite) matrix. The fact that the product
$\prod\limits_{i=1}^{n}U_{i}(X_{i})$ is unitary prove our claim.
We get that
$$\mu(X)=\mu(\prod\limits_{i=1}^{n}U_{i}(X_{i}))=\sum\limits_{i=1}^{n}\mu(U_{i}(X_{i}))$$
$$=\sum\limits_{i=1}^{n}\mu(X_{i}) \leq 4 \pi n.$$ Note that the
last equality follows from~\eqref{mu}. Finally, $H_{X}$ is
positive definite since it is a direct sum of positive definite
Hamiltonians (see formula~\eqref{prod. formula}).

It remains to prove the case $n=1$. We produce here a concrete
example. As argued above, we may assume that
$P=diag(\lambda,1/\lambda)$. Consider the function
$f(t)=\tan(\pi/4+at)$, where $|a|<\pi/4$ is chosen so that
$f(1)=\lambda$. We then have $$|f'|=|a(f^{2}+1)|<f^{2}+1.$$ Now,
let

$$U(t):= \left( \! \! \! \begin{array}{rr} \cos(2\pi t) & -\sin(2 \pi t) \\
\sin(2 \pi t)  &  \cos(2\pi t) \end{array} \right).  $$

\noindent Define the path $$X(t)=U(t)F(t)U(t)$$ where $F(t):=
diag(f(t),1/f(t))$. We claim that $\mu(X)=4 \pi$. Indeed,
$\mu(X)=\mu(U)+\mu(f)+\mu(U)=2 \pi +0+2 \pi =4 \pi$.

Finally we want to show that $X$ is positive.

\noindent Using~\eqref{prod. formula} we get:
$$H_{X}=H_{U}+UH_{F}U^{-1}+UF^{-1}H_{U}F^{-1}U^{-1},$$

\noindent and we know that $$H_{U}=2 \pi I; \text{  and  }H_{F}= \left( \begin{matrix} 0  & f'(t)/f(t) \\
f'(t)/f(t) & 0 \end{matrix} \right).$$ Thus $H_{X}$ is positive if
and only if $$U^{-1}H_{X}U=H_{U}+H_{F}+F^{-1}H_{U}F^{-1}$$

$$=2\pi (I+F^{-2})+H_{F}=\left( \begin{matrix} 2
\pi (1+f^{-2})  & f'(t)/f(t) \\
f'(t)/f(t) & 2 \pi (1+f^{2}) \end{matrix} \right)$$

is positive. Since the diagonal entries of the latter matrix are
positive, then the matrix is positive if and only if its
determinantis positive. Indeed we have:$$4
\pi^{2}(1+f^{2})(1+f^{-2})-(f'/f)^{2}=4
\pi^{2}(f^{2}+1/f)^{2}-(f'/f)^{2}>0$$ where in the last inequality
we have used the fact that $|f'|<1+f^{2}.$ This complete the proof
of the lemma.

\end{proof}

\subsubsection{Proof of Theorem~\ref{sp lemma} - Step 2}

 In what follows we denote by $C$ the constant which appears in
 the definition of the Maslov quasimorphism. That is:$$|\mu(fg)-\mu(f)-\mu(g)|\leq
 C$$for all $f,g \in \mathcal{S}$. It is not hard to prove the
 following fact (the proof appears in the literature see~\cite{Ba}).
\begin{lemma}\label{constants}
Let $\tilde{\mu}$ be the homogeneous quasimorphism of $\mu$.
 Then $$|\tilde{\mu}(X)-\mu(X)| \leq C\text{  for all X   }\in\mathcal{S}. $$
 We then have $|\tilde{\mu}(X)-\tilde{\mu}(Y)| \leq 4C$.
 \end{lemma} We will denote the constant $4C$ by $C_{1}$. The main
 lemma of step-2 is the following.
 \begin{lemma}\label{positive Y}Every element $Y \in \mathcal{S}$ for which
 $\mu({Y})\geq 6 \pi n+C$ is positive.
 \end{lemma}
 \begin{proof} Let $Y(t)=P(t)V(t)$ be the polar decomposition of
 $Y$. Let $X$ be the positive path, guaranteed by the main lemma of
 step-1, ending at $P(1).$ Define $$Z(t)=X^{-1}(t)P(t).$$ Note
 that $Z$ is a closed path, moreover we have \begin{equation}\label{eq.1}\mu(Z)\geq
 \mu(P)-\mu(X)-C \geq -4 \pi n-C\end{equation} where in the first inequality
 we have used the quasimorphism property, and in the second
 inequality we have used the fact that for every $X   \in\mathcal{S}
 $ we have $\mu(X)=-\mu(X^{-1}).$ From~\eqref{eq.1} we
 get:$$\mu(X^{-1}Y)=\mu(X^{-1}PV)=\mu(ZV)$$
 $$=\mu(Z)+\mu(V)=\mu(Z)+\mu(Y)>2 \pi n.$$ Now, since $X^{-1}Y(1)$
 is unitary, then $X^{-1}Y$ is homotopic to its unitary
 projection, we denote by $T$. So we have $\mu(T)=\mu(Y)>2 \pi n$
 and from~\ref{main linear lemma} it follows that $T\geq\Id$ in the
 unitary sense and thus in the symplectic sense. Consequently
 $X^{-1}Y$is positive and thus $Y=XX^{-1}Y$ is positive.

 \end{proof}

 We are now ready to prove the main theorem:

 We wish to compute $\gamma(X,Y)$. We claim that a sufficient
 condition for $X^{r}\geq Y^{s}$ to hold is
 \begin{equation}\label{condition}
 r\tilde{\mu}(X)-s\tilde{\mu}(Y)\geq 6 \pi
 n+2C+C_{1}.\end{equation}

 Indeed, if this condition holds, then $$\mu(X^{r}Y^{-s}) \geq
 \tilde{\mu}(X^{r}Y^{-s})-C \geq \tilde{\mu}(X^{r}) +
 \tilde{\mu}(Y^{-s})-C-C_{1}$$ $$=r
 \tilde{\mu}(X)-s\tilde{\mu}(Y)-C-C_{1} \geq 6 \pi n+C.$$

 Where we have used in the first and the second inequalities
 lemma~\ref{constants}, and in the third inequality we have used
 our assumption.

 Now, by lemma~\ref{positive Y}, we see that $X^{r}Y^{-s} \geq \Id$
 or that  $X^{r} \geq Y^{s}.$ Since the r.h.s.
 of~\eqref{condition} is independent of $r$ and $s$, then for
 every $\varepsilon \geq 0$ we can find large $r$ and $s$, which satisfy~\eqref{condition},
 such that
 $$0<\frac{r}{s}(\tilde{\mu}(X))-\tilde{\mu}(Y)<\varepsilon.$$
 Since $\gamma^{*}(X,Y) \leq r/s$, we deduce that
 $$\gamma^{*}(X,Y)=\gamma(X,Y) \leq \frac{\tilde{\mu}(Y)}
{\tilde{\mu}(X))}.$$ The fact that $\gamma(X,Y)\gamma(Y,X)\geq 1$
finishes the proof of theorem~\ref{sp lemma} and thus of
theorem~\ref{main linear theorem}.

\section{Proof of Theorems ~\ref{main theorem 1} and ~\ref{main theorem}}\label{quant
section}

Before we start proving the main theorems
 we add here one more
basic construction. That is, the construction of vector fields
generating the elements of the subgroup, of the contactomorphisms
group, $Quant$ (recall that $Quant$ is the subgroup of all
contactomorphisms which preserve the contact form, see
subsection~\eqref{pre. hofer}). The purpose of this construction
is to give an intuitive geometrical description of the group
$quant$.

submitted

Given  $(M,\omega)$ a symplectic manifold, such that $\omega\in
H^{2}(M,\mathbb{Z})$, and  $(P,\xi)$ its prequantization (as we
remarked in~\ref{pre. hofer} $P$ carries the structure of a
principal $S^{1}$ bundle over $M$) recall that in such a case the
manifold $P$ carries a 1-form, $\alpha$, globally defined on $P$,
such that $\xi=\{v\in T_{a}P|\alpha(v)=0,\text{ }\forall a\in
P\}.$ Moreover $\alpha$ is a connection form on $P$ such that
$\omega$ is its curvature . That is $d\alpha=p^{*}\omega$ where
$p:P\rightarrow M$ is the fiber bundle projection from $P$ to $M$
( See~\cite{Na} and~\cite{Sp} chapter 8 as a source on connections
and curvature of principal bundles). This means that for every
$a\in P$ we have the smooth decomposition
\begin{equation}\label{conn.decom.}T_{a}P=\xi_{|a}\oplus V_{a}\end{equation} where $V_{a}$ is the vertical
subspace of $T_{a}P$ defined canonically as the subspace tangent
to the fiber over the point $p(a)$, at the point $a.$ Equivalently
$V_{a}$ can be defined as the restriction to $a$ of all vector
fields which are the image under the Lie algebras homomorphisms
between the Lie algebra of $S^{1}$ into the Lie algebra of vector
fields of $P$ under the $S^{1}$ action (recall that this Lie
algebra homomorphism is induced by the action of the Lie group
$S^{1}$ on
$P$).

Now, it can be easily verified that for every $x \in M$ and $a\in
p^{-1}(x)$ $p_{*}(\xi_{|a})$ is an isomorphism between $T_{x}M$
and $\xi_{|a}$. Assume now that $p^{*}F$ is a contact Hamiltonian
on $P$ which is a lift of an Hamiltonian $F$ on $M$ (see
subsection~\ref{pre. hofer}).

Using decomposition~\eqref{conn.decom.} one can describe the
contact vector field obtained from $p^{*}F$ as follows: Denote by
$X_{F}(x)$ the Hamiltonian vector field obtained by $F$ at the
point $x$. Let $a \in p^{-1}(x)$ be a point above $x$. Then the
contact vector vector field $X_{p^{*}F}$ at the point $a$ is the
sum \begin{equation}X_{p^{*}F}=p_{*}^{-1}(X_{F}(x))\oplus
v\end{equation} where $p_{*}^{-1}(X_{F}(x))\in \xi_{|a}$ and $v$
is the unique vector in $V_{a}$ determined by the condition
$\alpha(v)=F(x).$ One can say that the horizontal component of the
contact vector field is determined by the (symplectic) Hamiltonian
vector field and its vertical componenet is the measure of
transversality (determined by $F$) to the horizontal field $\xi$.

Take for example any constant time dependant Hamiltonian
$c(t):M\rightarrow \mathbb{R}.$ Then in this case we have
$p_{*}^{-1}(X_{F}(x))=0$ and the vertical component of the contact
vector field is $v(t)=\alpha(c(t))$. It is easily verified that
the dynamics on $P$ after time $t$ is $$a\mapsto e^{ic(t)}(a)$$
where we have used the notation from subsection~\ref{pre. hofer}.

\subsection{Proof of Theorem~\ref{main theorem 1}}
We start with the following lemma in which we establish a
connection between the Hofer's metric and the partial order.


\begin{lemma} Let $f$ be the time-1-map of a flow generated by
Hamiltonian $F \in \mathcal{F}$. Let $\widetilde{F}$ be its lift
to the prequantization space and let $\tilde{f}$ be its
time-1-map. Then we have the formulas:
\begin{equation}\label{onee}\|f\|_{+}=\inf\{s|\text{ }e^{is}\geq\tilde{f}\}
\end{equation}

\begin{equation}\label{two}\|f\|_{-}=\inf\{s|\text{ }\tilde{f}\geq e^{-is}\}
\end{equation}
\end{lemma}

\begin{proof}
 We prove here formula~\eqref{onee}, formula~\eqref{two} is
proved along the same lines. Assume that
$$e^{is}\geq \tilde{f}\Leftrightarrow e^{is}\tilde{f}^{-1} \geq
\Id $$

This means that $\exists H_{1} \geq 0$ which generates the element
$e^{is}\tilde{f}^{-1}$.

\noindent Moreover we have that the contact Hamiltonian
$H=H_{1}-s$ connects $\Id$ to $\tilde{f}^{-1}$ which implies that
$-H$ connects $\Id$ to $\tilde{f}$. Now the fact that $H_{1} \geq
0 $ gives us that $-H \leq s$ $\Leftrightarrow$
 $\max(-H) \leq s$.

 \noindent Remembering that $-H$
generates $\tilde{f}$ (note that $H \in \mathcal{F}$, this can be
seen by using the Calabi-Weinstein invariant in a similar way to
the way we have use it in lemma~\ref{cw}) we have
$$\|f\|_{+} \leq \inf\{s|\text{
}e^{is}\geq\tilde{f}\}.$$

\noindent On the other hand assume that $s \geq \|f\|_{+}.$ This
means that $\exists F$ such that $F \in \mathcal{F}$ and $\max F
\leq s$. Define $-H=F.$ So we have $H+s \geq 0$ which implies that
$e^{is}\tilde{f}^{-1}\geq \Id$ which implies that $e^{is} \geq
\tilde{f}.$

\noindent Thus we conclude: $$\|f\|_{+} \geq \inf\{s|\text{
}e^{is}\geq\tilde{f}\}$$  which is the desired.

\end{proof}

As a result we have the following corollary.

\begin{corollary}\label{polt.lemma}Let $F:M\rightarrow \mathbb{R}$ be an Hamiltonian.
 Let $\{f_{t}\}$ be its flow and let $f$ be the time-1-map of the
flow. Let $\widetilde{F}$ be the lift of $F$ and denote by
$\tilde{f}$ its time-1-map. Then
\begin{equation}\label{1 of main} e^{is} \widetilde{f} \geq \Id
\Longleftrightarrow \|f\|_{-} \leq s
\end{equation} and \begin{equation}\label{2 of main} e^{is}
\widetilde{f} \leq \Id \Longleftrightarrow \|f\|_{+} \leq -s
\text{
  .}
\end{equation}
\end{corollary}
\begin{proof}
 We use formula~\eqref{onee} to derive formula~\eqref{2 of
 main}. The derivation of formula~\eqref{1 of main} from
 formula~\eqref{two} can be shown in the same way.

 For the first direction note that $$\|f\|_{+} \leq -k \Rightarrow
 \inf\{s|\text{ }e^{is}\geq\tilde{f}\} \leq -k$$

 $$\Rightarrow \text{    } \Id \geq e^{-i(-k)}\tilde{f} \text{   }
 \Rightarrow \Id \geq e^{ik} \tilde{f}.$$

 We now show the other direction. $$e^{ik}\tilde{f} \leq \Id \text{  }
 \Rightarrow \text{   } -k \geq \inf\{s|\text{
 }e^{is}\geq\tilde{f}\}$$
 $$\Rightarrow \|f\|_{+} \leq -k \text{,  as required  .}$$

 \end{proof}

We now prove the following formulas,
 ($t,s,\tilde{f}$ are
as in the theorem).

\begin{equation}\label{first}\gamma(e^{it},e^{is}\tilde{f})=\frac{s+\|f\|_{+,\infty}}{t}\end{equation}

\begin{equation}\label{second}\gamma(e^{is}\tilde{f},e^{it})=\frac{t}{s-\|f\|_{-,\infty}}.\end{equation}

We will prove formula~\eqref{first}. Formula~\eqref{second} is
proved along the same lines. First, recall that the function
$\gamma$ can be defined alternatively as in~\eqref{gamma*}, which
is the definition we use here. So assume that
$$c:=\frac{r}{p}\geq \gamma(e^{it},e^{is}\tilde{f})$$
$\Rightarrow$ $$e^{itr} \geq e^{isp}\tilde{f}^{p}$$ (here we use
the fact that all maps commutes, and the alternative definition of
$\gamma$) $$\Rightarrow e^{i(tr-sp)}\geq \tilde{f}^{p} \Rightarrow
tr-sp \geq \|f^{p}\|_{+}$$ (here we use formula~\eqref{onee})
$$\Rightarrow t \frac{r}{p} \geq \frac{\|f^{p}\|_{+}}{p}+s
\Rightarrow c \geq \frac{s+\frac{\|f^{p}\|_{+}}{p}}{t}$$

since $p$ can be chosen big as we want, we conclude that $$ c \geq
\frac{s+{\|f\|_{+,\infty}}}{t}$$ $\Rightarrow$
$$\gamma(e^{it},e^{is}\tilde{f}) \geq
\frac{s+\frac{\|f^{p}\|_{+}}{p}}{t}.$$

On the other hand assume that $$ c \geq
\frac{s+{\|f\|_{+,\infty}}}{t}.$$ Then there exists a sequence of
positive real numbers, $\{k_{n}\}$, such that $k_{n}\rightarrow
\infty$ and
$$c \geq \frac{s+\frac{\|f^{k_{n}}\|_{+}}{k_{n}}}{t}.$$ We
conclude that
$$tk_{n}c-sk_{n} \geq \|f^{k_{n}}\|_{+} \Rightarrow e^{i(\frac{trk_{n}}{p}-sk_{n})} \geq
\tilde{f}^{k_{n}}$$

\noindent $\Rightarrow$
\begin{equation}\label{frac}{e^{it}}^{\frac{rk_{n}}{p}} \geq
(e^{is}\tilde{f})^{k_{n}}\end{equation}

Now choose a sequence $\alpha_{n}$ such that $0 \leq \alpha_{n}
\leq 1$ and $\frac{rk_{n}}{p}+\alpha_{n} \in \mathbb{N}.$

From inequality~\eqref{frac} and the choice of $\alpha_{n}$ we
get:
$$e^{it({\frac{rk_{n}}{p}+\alpha_{n}})} \geq
(e^{is}\tilde{f})^{k_{n}}.$$

From the last inequality and the definition of $\gamma$ we use
here we get,

$$\frac{ \frac{rk_{n}}{p}+\alpha_{n}}{k_{n}} \geq
\gamma(e^{it},e^{is}\tilde{f})$$

 $\Rightarrow$

$$\frac{r}{p}+\frac{\alpha_{n}}{k_{n}} \geq
\gamma(e^{it},e^{is}\tilde{f}).$$

\noindent Now, since $\lim \limits_{n \rightarrow \infty}
\frac{\alpha_{n}}{k_{n}}=0$ we get that

$$c=\frac{r}{p} \geq \gamma(e^{it},e^{is}\tilde{f})$$ which is
what we need.

At this point we remark that due to~\eqref{first} and~\eqref
{second} and the fact that $\gamma(f,g)\gamma(g,f) \geq 1$ for
every $f$ and $g$, we infer that the r.h.s of~\eqref{log} is
defined.

Now we can actually calculate $K(e^{it},e^{is}\tilde{f})$. From
this calculation we will derive formula~\eqref{log}.

  By the very definition of $K$ and formulas~\eqref{first}
and~\eqref{second} we have
$$K(e^{it},e^{is}\tilde{f})=\max\{\log (s+\|f\|_{+,\infty})-\log
t\text{, }\log t- \log (s-\|f\|_{-,\infty})\}$$

$$=\frac{\log
(s+\|f\|_{+,\infty})- \log (s-\|f\|_{-,\infty})}{2}+\frac{|2\log t
-(\log (s+\|f\|_{+,\infty})+ \log (s-\|f\|_{-,\infty}))|}{2}$$ Now
clearly this expression attains its infimum when

$$t=e^{\frac{-(\log (s+\|f\|_{+,\infty})+ \log
(s-\|f\|_{-,\infty}))}{2}}.$$

Thus we conclude that the l.h.s of~\eqref{log} equals
$$=\frac{\log (s+\|f\|_{+,\infty})- \log
(s-\|f\|_{-,\infty})}{2}$$ or simply
$$\frac{1}{2}\log\frac{s+\|f\|_{+,\infty}}{s-\|f\|_{-,\infty}}$$
which is the desired.

\subsection{Proof of Theorem~\ref{main theorem}}

We begin with the following lemma.

\begin{lemma} Let $F \in \texttt{F}$ be an Hamiltonian. Then we
have the formulas:\begin{equation}\label{first eq.}\max
F=\|f\|_{+}\end{equation} and \begin{equation}\label{second
eq.}-\min F=\|f\|_{-}.\end{equation}
\end{lemma}

\begin{proof}

Our starting point is the following fact which is due to
Polterovich which can be easily deduced from~\cite{P-1}.

\textbf{Fact.} Let $(M,\omega)$ be a symplectic manifold. Let $L$
be a closed Lagrangian in $M$ with the stable Lagrangian
intersection property. Moreover, let $F$ be an autonomous
Hamiltonian defined on $M$ such that $F \in \mathcal{F}$ (see
subsection~\ref{pre. hofer}) and for some positive constant $C$ we
have $F_{|L} \geq C$. Denote by $f$ the time-1-map defined by $F$.
Then we have
\begin{equation}\label{p}\|f\|_{+} \geq C.\end{equation}

Now let $L_{\alpha}$ be a Lagrangian of the family of Lagrangians
which foliate $M_{0}$ as in theorem~\ref{main theorem}. Assume
that $F_{|L_{\alpha}} \geq C_{\alpha}.$ Then using~\eqref{p} we
have the following double inequality.
$$\max F\geq \|f\|_{+} \geq C_{\alpha}\text{, }\forall \alpha$$ where in the first
inequality we have used the very definition of the norm $\|\text{
}\|_{+}$. Now due to the fact that the set of Lagrangians foliate
an open dense set in $M$ we get that
$$\max F=\max \limits_{\alpha} \max \limits_{L_{\alpha}}F.$$ From this
we get we get formula~\eqref{first eq.}.

In the same manner we obtain formula \eqref{second eq.}.
\end{proof}


 \noindent We now calculate the relative growth on elements of $V$ (see~\ref{pre. hofer} for the definition of $V$).

 \noindent Let
$F,G \in \texttt{F}$. Denote by $\varphi_{F} \text{  the time-1-
map of }F$ and by $\varphi_{G} \text{  the time-1-map of }G.$ Let
$\widetilde{\varphi}_{F}$, $\widetilde{\varphi}_{G}$ be their lift
to the prequantization space. Now let $s,t$ be any real numbers
such that
$$e^{is}\widetilde{\varphi}_{F},e^{it}\widetilde{\varphi}_{G}\geq
\Id.$$ Note that this means that
$e^{is}\widetilde{\varphi}_{F},e^{it}\widetilde{\varphi}_{G}$ are
in $V$ and of course all elements of $V$ can be characterize in
this way. Now, by the very definition of the relative growth we
have
$$\gamma(e^{is}\widetilde{\varphi}_{F},e^{it}\widetilde{\varphi}_{G})
=\lim\limits_{n\rightarrow\infty}\frac{\gamma_{n}(e^{is}\widetilde{\varphi}_{F},e^{it}\widetilde{\varphi}_{G})}{n}.$$

We calculate $\gamma$ by a direct calculation of the functions
$\gamma_{n}$ on elements of $V$.

By the very definition of $\gamma_{n}$ we have

$$\gamma_{n}(e^{is}\widetilde{\varphi}_{F},e^{it}\widetilde{\varphi}_{G})=
\inf\{m|\text{ }e^{ims}{\widetilde{\varphi}}^{m}_{F}\geq
e^{int}{\widetilde{\varphi}}^{n}_{G}\}$$

$$=\inf\{m|\text{
}e^{i(ms-nt)}\widetilde{\varphi}_{mF-nG} \geq \Id\}=\inf\{m|\text{
}\|\varphi_{mF-nG}\|_{-}\leq ms-nt\}$$

(where in the second equality we have used the fact
  that the functions of the union $\texttt{F}\cup \{e^{is}\}_{s\in \mathbb{R}}$ are all Poisson
commutes, in the third equality we have used
corollary~\ref{polt.lemma})

$$=\inf\{m|\text{
}-\min(mF-nG)\leq ms-nt\}.$$ this follows  from
formula~\eqref{second eq.}.

 So we need $m$ that will satisfy
$$\max(nG-mF) \leq ms-nt\Leftrightarrow  $$

 $$ nG-mF\leq ms-nt \Leftrightarrow n(G+t)-m(F+s) \leq 0$$
$$\Leftrightarrow n(G+t) \leq m(F+s) $$

dividing both sides of the last inequality by the positive
function $F+s$ (recall that $F+s$ generates the dominant
$e^{is}\widetilde{\varphi}_{F}$) we get
 $$\Leftrightarrow \max (\frac{G+t}{F+s})n \leq m.$$ Thus we have

$$\gamma_{n}(e^{is}\widetilde{\varphi}_{F},e^{it}
 \widetilde{\varphi}_{G})-1 \leq \max(\frac{G+t}{F+s})n \leq \gamma_{n}(e^{is}\widetilde{\varphi}_{F},e^{it}
 \widetilde{\varphi}_{G})$$ thus

 $$\frac{\gamma_{n}(e^{is}\widetilde{\varphi}_{F},e^{it}
 \widetilde{\varphi}_{G})-1}{n} \leq \max(\frac{G+t}{F+s}) \leq \frac{\gamma_{n}(e^{is}\widetilde{\varphi}_{F},e^{it}
 \widetilde{\varphi}_{G})}{n}.$$

 So we have the formula

 \begin{equation}\label{gamma for.}\gamma(e^{is}\widetilde{\varphi}_{F},e^{it}\widetilde{\varphi}_{G})
 =\max(\frac{G+t}{F+s}).
 \end{equation}

 Now, let $\widetilde{f},\widetilde{g}\in V$ generated by the
 Hamiltonians $\widetilde{F},\widetilde{G}$ respectively. Then
 according to formula~\eqref{gamma for.} we have
 \begin{equation}\label{K}
 K(\widetilde{f},\widetilde{g})=\max|\log\widetilde{F}-\log\widetilde{G}|
\end{equation}

Using formula~\eqref{K}, for $K$, we define the isometry of
$\{\texttt{F},\|\text{ }\|_{max}\}$ to $Z.$

Let $\widetilde{f} \in V $ generated by the Hamiltonian
$\widetilde{F}$. Then the correspondence
$\widetilde{f}\leftrightarrow \log\widetilde{F}$ clearly induces
the required isometric imbedding of $\{\texttt{F},\|\text{
}\|_{max}\}$ into $Z$. This conclude the proof of this part of the
theorem.

\subsubsection{Examples}

\textbf{Example 1.} Consider the $2n$ dimensional standard
symplectic torus. That is $T^{2n}=\mathbb{R}^{2n}/\mathbb{Z}^{2n}$
with the symplectic structure $dP \wedge dQ=\sum
\limits_{i=1}^{n}dp_{i} \wedge dq_{i}.$ We now foliate the tours
by a family of Lagrangians which depends only on the $P$
coordinate. This way we exhibit the tours as a Lagrangian
fiberation parameterize by an $n$ tours (the $P$ coordinate which
in the notations of theorem~\ref{main theorem} is the parameter
$\Lambda$) and with Lagrangian fibers (the $Q$ coordinate). It is
known that this family of Lagrangian has the stable Lagrangian
intersection property. Now, define the family of autonomous
functions on $T^{2n}$ which depends only on the $P$ coordinate.
Thus all the conditions of theorem~\ref{main theorem} are
satisfied (we remained the reader that all functions of this
family are Poisson commuting). \noindent We conclude that the
metric space $Z$ of a prequantization space of the standard $2n$
dimensional symplectic torus contains an infinite dimensional
metric space. Denote by $C^{\infty}_{nor}(T^{n})$ the space of
smooth normalized functions on the $n$-dimensional torus. Then we
have the isometric injection:$$(C^{\infty}_{nor}(T^{n}),\|\text{
}\|_{\max})\hookrightarrow Z.$$

\begin{figure}[h!t]
\begin{minipage}{4.5cm}
\begin{center}
\epsfig{file=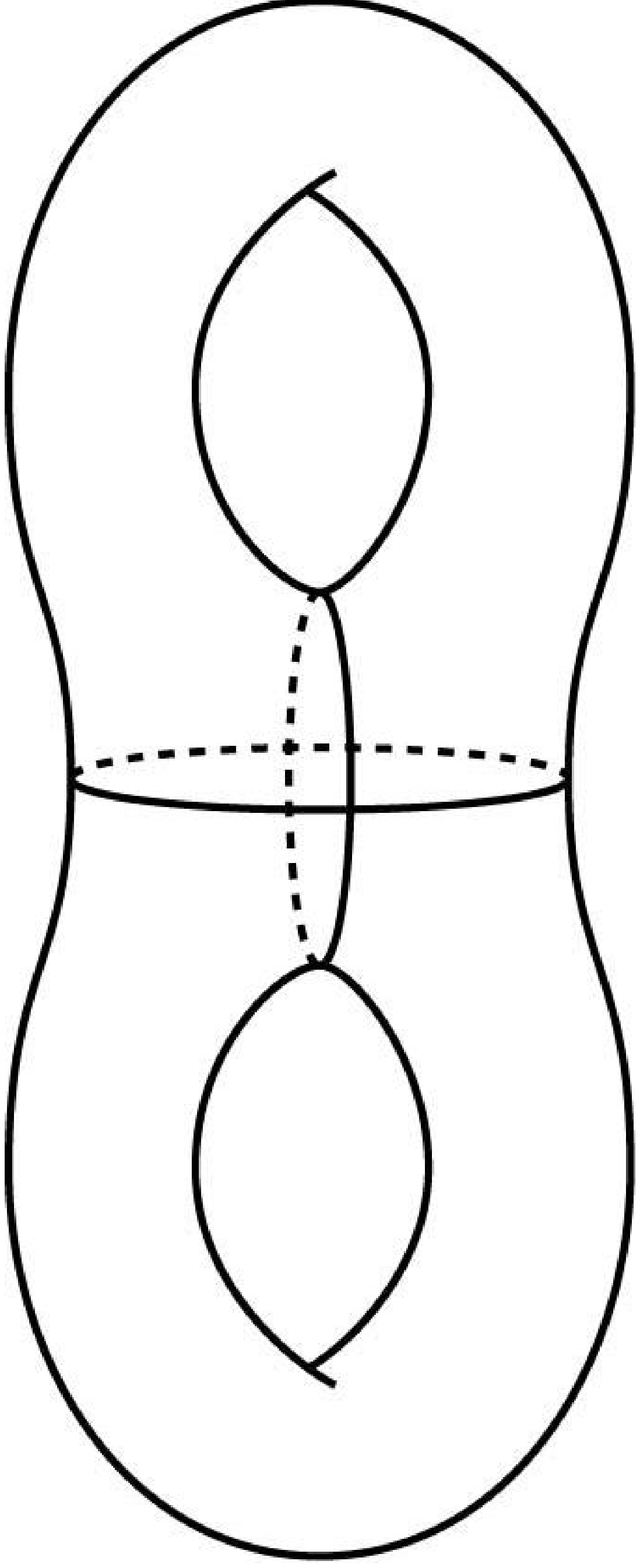,width=1.9cm} \caption{The curves we
remove from $\Sigma$. } \label{surface_2}
\end{center}
\end{minipage}
\hfill
\begin{minipage}{4.5cm}
\begin{center}
\epsfig{file=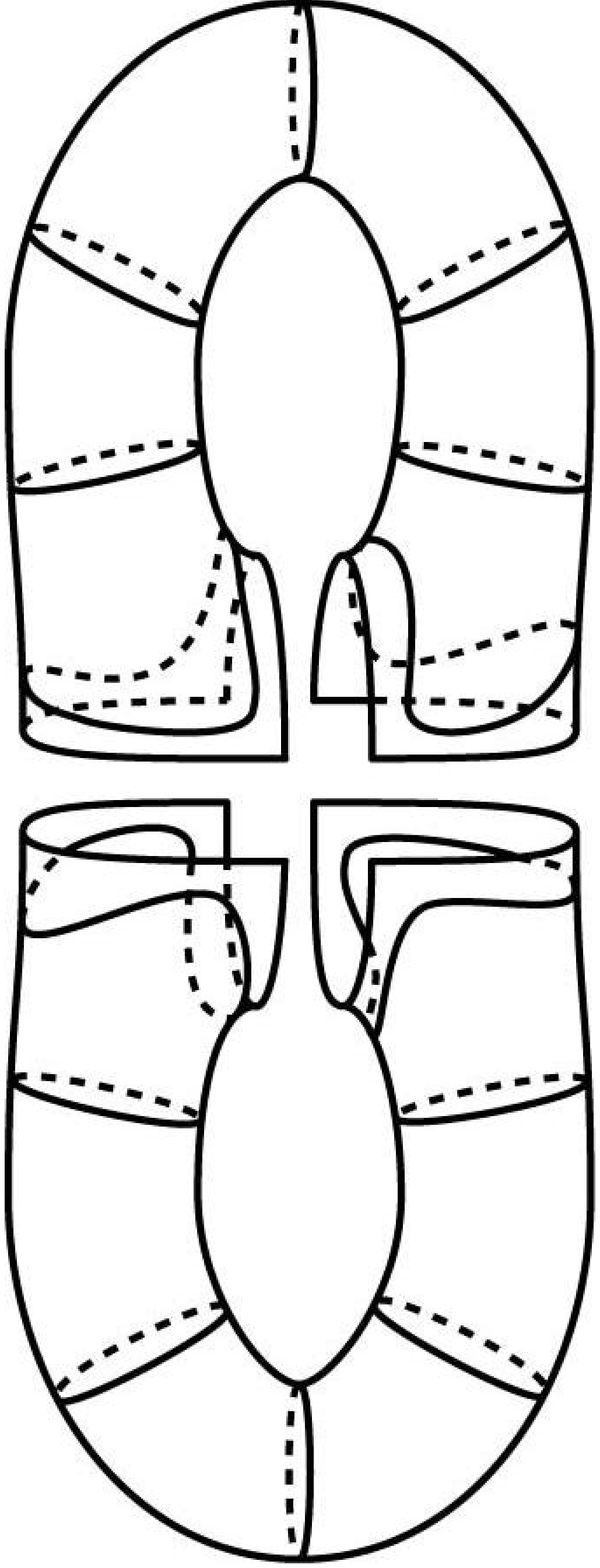,width=1.9cm} \caption{The foliation of
$\Sigma_{0}$. } \label{foliation}
\end{center}
\end{minipage}
\end{figure}

\bigskip
\textbf{Example 2.} Let $\Sigma$ be any surface of genus greater
 or equal 2. It is known that non contractible loops on $\Sigma$
 has the stable Lagrangian property (see for example~\cite{P-S}). For such surfaces one can
 foliate a subset, $\Sigma_{0}$, of $\Sigma$ by a family of disjoint non
 contractible closed loops such that $\Sigma_{0}$ is an open dense
 subset of $\Sigma$ (actually $\Sigma\setminus \Sigma_{0}$
 is a finite collection of closed arcs-see the figures 1 and 2).
 Note that all the conditions of theorem~\eqref{main theorem} are satisfied. In the figures 1 and 2 we show an
 example of such family of Lagrangians for surface of genus equal
 2.

\bigskip

\textbf{Acknowledgements.} I want to thank deeply  Leonid
Polterovich for reaching an helping hand each time it was needed,
 in all aspects of academic research activity: support, ideas, and
encouragement. I also want to thank Yasha Eliashberg for inviting
me to Stanford university and for sharing with me his opinion on
this work, and his kind hospitality. The "linear" part of this
work has started at the beginning of my Ph.D studies and has
benefit a lot from many discussions with Assaf Goldberger. I thank
him very much for his help.

\end{document}